\documentclass[12pt]{article}
\usepackage[leqno]{amsmath}
\usepackage{amssymb,mathrsfs,bm}
\usepackage{hyperref}
\usepackage{mathtools}
\usepackage{fullpage}
\usepackage{epic}
\usepackage{eepic}
\usepackage[english]{babel}
\usepackage[all]{xy}
\usepackage{enumerate}
\usepackage[T1]{fontenc}
\usepackage[latin1]{inputenc}
\usepackage{color}
\usepackage{soul}
\usepackage{mathdots}
\usepackage{mathabx}
\usepackage{dsfont}

\newtheorem{theo}{Theorem}
\newtheorem{lem}{Lemma}

\newtheorem{prop}{Proposition}

\DeclareMathOperator{\Tra}{Tr}

\DeclareMathOperator{\vol}{vol}

\DeclareMathOperator{\cF}{\mathcal{F}}

\DeclareMathOperator{\cA}{\mathcal{A}}

\DeclareMathOperator{\gl}{\mathfrak{gl}}

\DeclareMathOperator{\fh}{\mathfrak{h}}

\DeclareMathOperator{\cO}{\mathcal{O}}

\DeclareMathOperator{\rs}{rs}

\DeclareMathOperator{\GL}{\mathrm{GL}}

\DeclareMathOperator{\val}{val}
\DeclareMathOperator{\Mat}{\mathrm{Mat}}
\DeclareMathOperator{\bun}{\mathbf{1}_{\fh_n(\cO_F)}}
\DeclareMathOperator{\bgln}{\mathbf{1}_{\gl_n(\cO_F)}}
\DeclareMathOperator{\rU}{\mathrm{U}}
 
\newcommand{\eq}[1][r]
   {\ar@<-3pt>@{-}[#1]
    \ar@<-1pt>@{}[#1]|<{}="gauche"
    \ar@<+0pt>@{}[#1]|-{}="milieu"
    \ar@<+1pt>@{}[#1]|>{}="droite"
    \ar@/^2pt/@{-}"gauche";"milieu"
    \ar@/_2pt/@{-}"milieu";"droite"}

\title{A new proof of Jacquet-Rallis's fundamental lemma}
\author{Rapha\"el Beuzart-Plessis \protect\footnote{The project leading to this publication has received funding from Excellence Initiative of Aix-Marseille University-A*MIDEX, a French ``Investissements d'Avenir" programme.}}

\begin{document}

\maketitle

\begin{abstract}
We give a new proof of the so-called Lie algebra version of Jacquet-Rallis's fundamental lemma for local non-Archimedean fields of characteristic zero. This proof is local and based on a previous result of W. Zhang on the compatibility of smooth transfer with a (partial) Fourier transform.
\end{abstract}

\tableofcontents

\section{The linearized Jacquet-Rallis's fundamental lemma: statement of the main result}

Let $E/F$ be an unramified quadratic extension of local non-Archimedean fields of characteristic zero. We denote by $N$ and $\Tra_{E/F}$ the norm and trace of this extension respectively. Let $\sigma$ be the nontrivial $F$-automorphism of $E$ and $\cO_F$, $\cO_E$ the rings of integers of $F$ and $E$ respectively. For a matrix $M$, we denote by ${}^t M$ its transpose and, if the entries of $M$ are valued in $E$, by $M^\sigma$ the matrix obtained by applying $\sigma$ to each entry. Let $n\geqslant 1$ be an integer. All algebraic varieties will be tacitly assumed to be defined over $F$ and, by abuse of notation, we denote by $\gl_{n,E}$ and $\GL_{n-1,E}$ the Weil's restriction of scalars from $E$ to $F$ of $\gl_n$ and $\GL_{n-1}$ (viewed as varieties over $E$) respectively. We set
$$\displaystyle \fh_n:=\left\{X\in \gl_{n,E}\mid X={}^t X^\sigma  \right\},$$
$$\displaystyle \rU_{n-1}:=\left\{g\in \GL_{n-1,E}\mid {}^t g^\sigma g=I_{n-1} \right\}.$$
We let $\rU_{n-1}$, resp. $\GL_{n-1}$, act on $\fh_n$, resp. $\gl_n$, by conjugation through the embedding $\rU_{n-1}\hookrightarrow \rU_n$, resp. $\GL_{n-1}\hookrightarrow \GL_n$, given by $g\mapsto \begin{pmatrix} g & \\ & 1 \end{pmatrix}$. We also write, again by abuse of notation, $\rU_{n-1}(\cO_F)$ and $\fh_n(\cO_F)$ for $\rU_{n-1}(F)\cap \GL_{n-1}(\cO_E)$ and $\fh_n(F)\cap \gl_n(\cO_E)$ respectively.

An element $X\in \fh_n$ (resp. $Y\in \gl_n$) is said to be {\em relatively regular semi-simple} if its stabilizer in $\rU_{n-1}$ (resp. in $\GL_{n-1}$) is trivial and its orbit is closed. We denote by $\fh_n^{\rs}\subset \fh_n$ and $\gl_n^{\rs}\subset \gl_n$ the open subschemes of relatively regular semi-simple elements. We say that two elements $X\in \fh_n^{\rs}(F)$ and $Y\in \gl_n^{\rs}(F)$ {\em match}, and we will write $X\leftrightarrow Y$, if they are conjugate under $\GL_{n-1}(E)$ inside $\gl_n(E)$.

Set
$$\displaystyle e_n=\begin{pmatrix} 0 \\ \vdots \\ 0 \\ 1 \end{pmatrix}\in \Mat_{n,1}(F),\; e_n^*=(0,\ldots,0,1)\in \Mat_{1,n}(F).$$
For every $X\in \gl_n^{\rs}(F)$, we define
$$\displaystyle v(X):=\val(\det(e_n^*,e_n^*X,\ldots,e_n^*X^{n-1})),$$
$$\displaystyle \omega(X):=(-1)^{v(X)},$$
where $\val$ stands for the normalized valuation on $F$ and the determinant above is nonzero because $X$ is relatively regular semi-simple.

For $X\in \fh_n^{\rs}(F)$, $Y\in \gl_n^{\rs}(F)$ and functions $f\in C_c^\infty(\fh_n(F))$, $f'\in C_c^\infty(\gl_n(F))$, following Jacquet and Rallis we introduce the (relative and twisted) orbital integrals
$$\displaystyle O(X,f)=\int_{\rU_{n-1}(F)} f(gXg^{-1})dg$$
and
$$\displaystyle O(Y,f')=\int_{\GL_{n-1}(F)} f'(gYg^{-1})\omega(gYg^{-1})dg=\omega(Y) \int_{\GL_{n-1}(F)} f'(gYg^{-1})(-1)^{\val(\det g)}dg$$
where the Haar measures are normalized by requiring that $\vol(\rU_{n-1}(\cO_F))=\vol(\GL_{n-1}(\cO_F))=1$. Note however that the orbital integral $O(Y,f')$ differs from the one introduced in \cite{JR} by the ``transfer factor'' $\omega(Y)$. It is actually more convenient for us to work with this slightly modified definition as it makes the function $Y\in \gl_n^{\rs}(F)\mapsto O(Y,f')$ obviously invariant by $\GL_{n-1}(F)$-conjugation.

We say that two functions $f\in C_c^\infty(\fh_n(F))$ and $f'\in C_c^\infty(\gl_n(F))$ {\em match}, or that they are transfer of each other, if for every pair $(X,Y)\in \fh_n^{\rs}(F)\times \gl_n^{\rs}(F)$ of matching elements we have
$$\displaystyle O(X,f)=O(Y,f'),$$
and for every element $Y\in \gl_n^{\rs}(F)$ not matching any element of $\fh_n^{\rs}(F)$ we have
$$\displaystyle O(Y,f')=0.$$

We denote by $\bun$ and $\bgln$ the characteristic functions of $\fh_n(\cO_F)$ and $\gl_n(\cO_F)$ respectively.

The goal of this note is to prove the following result (also known as Jacquet-Rallis's fundamental lemma for Lie algebras see \cite[Conjecture 2]{JR}).

\begin{theo}\label{theo 1}
The functions $\bun$ and $\bgln$ match. That is to say (repeating the above definition):
\begin{enumerate}[(i)]
\item For every pair $(X,Y)\in \fh_n^{\rs}(F)\times \gl_n^{\rs}(F)$ of matching elements, we have
$$\displaystyle O(X,\bun)=O(Y,\bgln).$$

\item For every $Y\in \gl_n^{\rs}(F)$ not matching any element of $\fh_n^{\rs}(F)$, we have
$$\displaystyle O(Y,\bgln)=0.$$
\end{enumerate}

\end{theo}

By \cite[\S 2.6]{Yu} this version of Jacquet-Rallis's fundamental lemma implies a group version which was subsequently used by W. Zhang \cite{Zha1} to prove the global Gan-Gross-Prasad conjecture for unitary groups under some local assumptions. The above theorem was established by Yun in \cite{Yu} for local fields of positive characteristic greater than $n$ (except for (ii) which was proved by Yun in an elementary and direct way without any assumption on the characteristic). Moreover, in the appendix to {\it loc. cit.}, J. Gordon has shown that the transfer principle of Cluckers and Loeser, which relies on model-theoretic techniques, applies to this case thus deducing the Jacquet-Rallis's fundamental lemma for characteristic zero local fields in sufficiently large residual characteristic. Therefore the above theorem is not new, except in small residual characteristic, but the proof that we give is substantially different from \cite{Yu} which was based on the same geometric and cohomological ideas as Ng\^o's proof of Langlands-Shelstad's fundamental lemma. Instead, we build our argument on the compatibility of the Jacquet-Rallis's transfer with a certain partial Fourier transform that was established by W. Zhang in \cite{Zha1} that we reformulate in terms of Weil's representations in Section \ref{Section Weil}. In particular, our argument is purely local and only based on harmonic-analytic tools (such as the Fourier transform).

Finally, we remark that by the recent work of Jingwei Xiao \cite{Xi} the Jacquet-Rallis fundamental lemma implies the (Lie algebra version of) endoscopic fundamental lemma of Langlands-Shelstad for unitary groups. Therefore, putting together this paper with \cite{Zha1} and \cite{Xi} gives an elementary proof of the fundamental lemma for unitary groups! Actually, in \cite{Xi} it is also shown that Zhang's results \cite{Zha1} directly imply that the usual endoscopic transfer for (Lie algebras of) unitary groups ``commutes'' (up to an explicit constant) with the Fourier transform. This compatibility to the Fourier transform was previously proved by Waldspurger \cite{Wa1},\cite{Wa2} in the general setting of endoscopy under the assumption of the fundamental lemma but a very nice result of Kazhdan and Varshavsky \cite{KV} states that actually the converse holds: compatibility of the transfer with the Fourier transform implies the endoscopic fundamental lemma. Therefore, the combination of the papers \cite{Zha1}, \cite{Xi} and \cite{KV} also give an alternative, purely local, proof of the endoscopic fundamental lemma for unitary groups although the arguments in \cite{KV} are arguably slightly less elementary than the one in this paper (in particular, being based on Deligne-Lusztig's theory).

\section{Induction and a first case of Theorem \ref{theo 1}}

We will prove Theorem \ref{theo 1} by induction on $n$, the case $n=1$ being trivial. Thus, we assume from now on that $n\geqslant 2$ and that the result is known for $n$ replaced by $n-1$.

Every element $X\in \gl_n(E)$ can be written in blocks
$$\displaystyle X=\begin{pmatrix} X' & b \\ c & \lambda \end{pmatrix}$$
where $X'\in \gl_{n-1}(E)$, $b\in\Mat_{n-1,1}(E)$, $c\in \Mat_{1,n-1}(E)$ and $\lambda\in E$. We set
$$\displaystyle q(X)=cb\in E.$$
It is straightforward to see that $q$ restricts to a regular $\rU_{n-1}$-invariant (resp. $\GL_{n-1}$-invariant) function on $\fh_n$ (resp. on $\gl_n$). Moreover, for matching elements $X\in \fh_n^{\rs}(F)$ and $Y\in \gl_n^{\rs}(F)$ we have $q(X)=q(Y)$. The induction hypothesis will be used in the proof of the next lemma.

\begin{lem}\label{lem 1}
\begin{enumerate}[(i)]
\item Let $X\in \fh_n^{\rs}(F)$ and $Y\in \gl_n^{\rs}(F)$ be two matching elements and assume that $\lvert q(X)\rvert=\lvert q(Y)\rvert\geqslant 1$. Then, we have
$$\displaystyle O(X,\bun)=O(Y,\bgln).$$

\item Let $Y\in \gl_n^{\rs}(F)$. Assume that $Y$ does not match any element of $\fh_n^{\rs}(F)$ and that $\lvert q(Y)\rvert\geqslant 1$. Then, we have
$$\displaystyle O(Y,\bgln)=0.$$
\end{enumerate}

\end{lem}

\noindent\ul{Proof}: We prove (i), the proof of (ii) being similar. When $\lvert q(X)\rvert=\lvert q(Y)\rvert> 1$, it is easy to see that $O(X,\bun)=O(Y,\bgln)=0$. Therefore we may assume that $q(X)=q(Y)=\mu\in \cO_F^\times$. Since the extension $E/F$ is unramified, there exists $\nu\in \cO_E^\times$ such that $\mu=N(\nu)$. Then, it is easy to see that up to conjugation by $\rU_{n-1}(F)$ and $\GL_{n-1}(F)$ respectively, $X$ and $Y$ have the following forms
$$\displaystyle X=\begin{pmatrix} X' & \nu e_{n-1} \\ \sigma(\nu) e_{n-1}^* & \lambda \end{pmatrix},\; Y=\begin{pmatrix} Y' & N(\nu)e_{n-1} \\ e_{n-1}^* & \lambda' \end{pmatrix}$$
where $X'\in \fh_{n-1}(F)$, $Y'\in \gl_{n-1}(F)$ and $\lambda,\lambda'\in F$. Since $X$ and $Y$ match, we have $\lambda=\lambda'$ and moreover we readily check that $X'$ and $Y'$ also match. By the induction hypothesis we have
\begin{align}\label{eq 1}
\displaystyle O(X',\mathbf{1}_{\fh_{n-1}(\cO_F)})=O(Y',\mathbf{1}_{\gl_{n-1}(\cO_F)}).
\end{align}
Therefore, the lemma would follow if we show that
\begin{align}\label{eq 2}
\displaystyle O(X,\bun)=\mathbf{1}_{\cO_F}(\lambda)O(X',\mathbf{1}_{\fh_{n-1}(\cO_F)}),
\end{align}
and
\begin{align}\label{eq 3}
\displaystyle O(Y,\bgln)=\mathbf{1}_{\cO_F}(\lambda)O(Y',\mathbf{1}_{\gl_{n-1}(\cO_F)}).
\end{align}
We will show \eqref{eq 2}, the proof of \eqref{eq 3} being similar. By definition, we have
$$\displaystyle O(X,\bun)=\int_{\rU_{n-1}(F)} \bun\begin{pmatrix}gX'g^{-1} & \nu ge_{n-1} \\ \sigma(\nu) e_{n-1}^* g^{-1} & \lambda \end{pmatrix}dg.$$
Let $g\in \rU_{n-1}(F)$ be such that $\bun\begin{pmatrix}gX'g^{-1} & \nu ge_{n-1} \\ \sigma(\nu) e_{n-1}^* g^{-1} & \lambda \end{pmatrix}\neq 0$. Then, $ge_{n-1}\in \Mat_{n-1,1}(\cO_E)$ and ${}^t(ge_{n-1})^\sigma ge_{n-1}=1$. The group $\rU_{n-1}(\cO_F)$ acts transitively on the set of vectors $v\in \Mat_{n-1,1}(\cO_E)$ such that ${}^t v^\sigma v=1$. Since the stabilizer of $e_{n-1}$ in $\rU_{n-1}$ is $\rU_{n-2}$ it follows that $g\in \rU_{n-1}(\cO_F)\rU_{n-2}(F)$ i.e. the above integral is supported in $\rU_{n-1}(\cO_F)\rU_{n-2}(F)$. By the choice of Haar measures on $\rU_{n-1}(F)$ and $\rU_{n-2}(F)$ and since the function $\bun$ is invariant by conjugation by $\rU_{n-1}(\cO_F)$, it follows that
\[\begin{aligned}
\displaystyle & O(X,\bun)=\int_{\rU_{n-2}(F)} \bun\begin{pmatrix}gX'g^{-1} & \nu e_{n-1} \\ \sigma(\nu) e_{n-1}^* & \lambda \end{pmatrix}dg \\
 & =\mathbf{1}_{\cO_F}(\lambda)\int_{\rU_{n-2}(F)} \mathbf{1}_{\fh_{n-1}(\cO_F)}(gX'g^{-1})dg=\mathbf{1}_{\cO_F}(\lambda)O(Y',\mathbf{1}_{\gl_{n-1}(\cO_F)}).
\end{aligned}\]
$\blacksquare$

\section{Weil representations and transfer}\label{Section Weil}

Let $\psi$ be a nontrivial unramified additive character of $F$ (i.e. whose conductor is $\cO_F$) and set
$$\displaystyle \psi_E=\psi\circ \Tra_{E/F}$$
(an additive character of $E$). We define a ``partial" Fourier transform $\cF$ on $C_c^\infty(\fh_n(F))$ (respectively on  $C_c^\infty(\gl_n(F))$) by
$$\displaystyle (\cF f)\begin{pmatrix} X' & b \\ {}^t b^\sigma & \lambda \end{pmatrix}=\int_{\Mat_{n-1,1}(E)}f\begin{pmatrix} X' & c \\ {}^t c^\sigma & \lambda \end{pmatrix} \psi_E({}^t c^\sigma b)dc$$
for $f\in C_c^\infty(\fh_n(F))$ and $\begin{pmatrix} X' & b \\ {}^t b^\sigma & \lambda \end{pmatrix}\in \fh_n(F)$ (respectively by
$$\displaystyle (\cF f)\begin{pmatrix} X' & b \\ c & \lambda \end{pmatrix}=\int_{\Mat_{n-1,1}(F)\times \Mat_{1,n-1}(F)}f\begin{pmatrix} X' & b' \\ c' & \lambda \end{pmatrix} \psi(c'b+cb')db'dc',$$
for $f\in C_c^\infty(\gl_n(F))$ and $\begin{pmatrix} X' & b \\ c & \lambda \end{pmatrix}\in \gl_n(F)$) where the Haar measures are autodual i.e. such that $\cF$ has order $4$. Notice that, with this normalization and because $\psi$ is unramified, we have $\cF \bun=\bun$ and $\cF \bgln=\bgln$.

We have isomorphisms
$$\displaystyle (\gl_{n-1}(F)\oplus F)\times (F^{n-1}\oplus F^{n-1})\simeq \gl_n(F),$$
$$\displaystyle (X',\lambda,b,c)\mapsto \begin{pmatrix} X' & b \\ {}^t c & \lambda \end{pmatrix}$$
and
$$\displaystyle (\fh_ {n-1}(F)\oplus F)\times E^{n-1}\simeq \fh_n(F),$$
$$\displaystyle (X',\lambda,b)\mapsto \begin{pmatrix} X' & b \\ {}^t b^\sigma & \lambda \end{pmatrix}$$
which induce identifications
$$\displaystyle C_c^\infty(\gl_n(F))=C_c^\infty(\gl_{n-1}(F)\oplus F)\otimes C_c^\infty(F^{n-1}\oplus F^{n-1})$$
and
$$\displaystyle C_c^\infty(\fh_n(F))=C_c^\infty(\fh_{n-1}(F)\oplus F)\otimes C_c^\infty(E^{n-1}).$$
Let $W$ denote the Weil representations of $SL_2(F)$ on $C_c^\infty(F^{n-1}\oplus F^{n-1})$ and $C_c^\infty(E^{n-1})$ associated to the restrictions of the quadratic form $q$ (see e.g. \cite[Theorem 4.8.5]{Bump}) that we extend to $C_c^\infty(\gl_n(F))$ and $C_c^\infty(\fh_n(F))$ by letting $SL_2(F)$ act trivially on $C_c^\infty(\gl_{n-1}(F)\oplus F)$ and $C_c^\infty(\fh_{n-1}(F)\oplus F)$. Since $(F^{n-1}\oplus F^{n-1},q)$ can be written as an orthogonal sum of hyperbolic planes and $(E^{n-1},q)$ as an orthogonal direct sum of quadratic planes of the form $(E,N_{E/F})$ and $E/F$ is unramified, we have (the point being that the Weil's constant $\gamma_\psi(q)$ equals one in both cases, this follows e.g. from \cite[Lemma 1.2 (i) and (iv)]{JL})
\begin{align}\label{eqA}
\displaystyle W\begin{pmatrix} 1 & t \\ 0 & 1 \end{pmatrix}f(X)=\psi(tq(X)) f(X),
\end{align}
\begin{align}\label{eqB}
\displaystyle W\begin{pmatrix}  0 & -1 \\ 1 & 0 \end{pmatrix}f=\cF f
\end{align}
for every $f\in C_c^\infty(\gl_n(F))$ (respectively $C_c^\infty(\fh_n(F))$), $t\in F$ and $X\in \gl_n(F)$ (respectively $\fh_n(F)$).

Let $\cA$ be the GIT quotient $\gl_n//\GL_{n-1}$. Then, there is a natural isomorphism $\cA\simeq \fh_n//\rU_{n-1}$ (see e.g. \cite[Proposition 2.2.2.1]{Chau}). Note that $q$ descends to a regular function on $\cA$ that we shall denote by the same letter. Moreover, the image of $\gl_n^{\rs}$ in $\cA$ is an open subscheme $\cA^{\rs}\subseteq \cA$ and the inverse image of any $a\in \cA^{\rs}(F)$ in $\gl_n(F)$ consists of one relatively regular semi-simple $\GL_{n-1}(F)$-orbit (\cite[Lemmas 2.1.5.1 \& 2.1.5.3]{Chau}). Similarly, the inverse image of any $a\in \cA^{\rs}(F)$ in $\fh_n(F)$ is either empty or consists of one relatively regular semi-simple $\rU_{n-1}(F)$-orbit (\cite[Proposition 2.2.4.1]{Chau}). Denote by $p_{\GL}: \gl_n(F)\to \cA(F)$ and $p_U:\fh_n(F)\to \cA(F)$ the natural maps. For $f\in C_c^\infty(\fh_n(F))$, $f'\in C_c^\infty(\gl_n(F))$ and every $a\in \cA^{\rs}(F)$ we set
$$\displaystyle O(a,f)=\left\{
    \begin{array}{ll}
        O(X_a,f) \mbox{ if } p_U^{-1}(a)\neq \emptyset \mbox{ and where } X_a\in p_U^{-1}(a), \\
        0 \mbox{ otherwise,}
    \end{array}
\right.$$
and
$$\displaystyle O(a,f')=O(Y_a,f') \mbox{ where } Y_a\in p_{\GL}^{-1}(a).$$
(Note that the orbital integrals $O(.,f)$, $O(.,f')$ being invariant by conjugation by $\rU_{n-1}(F)$ (resp. $\GL_{n-1}(F)$) the above definitions are indeed independent of the choices of $X_a$ and $Y_a$.)

We also note that two functions $f\in C_c^\infty(\gl_n(F))$ and $f'\in C_c^\infty(\fh_n(F))$ match if and only if
$$\displaystyle O(a,f)=O(a,f')$$
for every $a\in \cA^{\rs}(F)$.

Set
$$\displaystyle Orb(\gl_n)=\{a\in \cA^{\rs}(F)\mapsto O(a,f')\mid f'\in C_c^\infty(\gl_n(F)) \}$$
and
$$\displaystyle Orb(\fh_n)=\{a\in \cA^{\rs}(F)\mapsto O(a,f)\mid f\in C_c^\infty(\fh_n(F)) \}.$$
Then, both $Orb(\gl_n)$ and $Orb(\fh_n)$ are subspaces of $C^\infty(\cA^{\rs}(F))$ the space of all locally constant complex-valued functions on $\cA^{\rs}(F)$ (this follows from the fact that the restrictions of $p_U$ and $p_{\GL}$ to the regular semi-simple locus $\cA^{\rs}(F)$ are $F$-analytic locally trivial fibrations see e.g. \cite[Lemma 3.12]{Zha1}).

\begin{prop}
The Weil representations $W$ on $C_c^\infty(\gl_n(F))$ and $C_c^\infty(\fh_n(F))$ descend to representations of $SL_2(F)$ on $Orb(\gl_n)$ and $Orb(\fh_n)$ which coincide on the intersection.
\end{prop}
\noindent\ul{Proof}: As $SL_2(F)$ is generated by $\begin{pmatrix} 1 & t \\ 0 & 1 \end{pmatrix}$, $t\in F$, and $\begin{pmatrix}  0 & -1 \\ 1 & 0 \end{pmatrix}$ it suffices to check that the actions of these elements descent to the spaces of orbital integrals and coincide on their intersection. For elements of the form $\begin{pmatrix} 1 & t \\ 0 & 1 \end{pmatrix}$, $t\in F$, this follows directly from \eqref{eqA} and moreover we have
\begin{align}\label{eqC}
\displaystyle W\begin{pmatrix} 1 & t \\ 0 & 1 \end{pmatrix}\Phi(a)=\psi(tq(a)) \Phi(a)
\end{align}
for every $\Phi\in Orb(\gl_n)+Orb(\fh_n)$, $t\in F$ and $a\in \cA^{\rs}(F)$. For $\begin{pmatrix}  0 & -1 \\ 1 & 0  \end{pmatrix}$, this follows from \eqref{eqB} and the following result of Zhang (see \cite[Theorem 4.17]{Zha1}, where it is proved that the transfer ``commutes" with $\cF$ up to a constant $\nu$, and \cite[\S 3.4.2]{Chau} for the precise computation of $\nu$ which is equal to $1$ in the unramified case we are considering):
\begin{theo}[Zhang]
If two functions $f'\in C_c^\infty(\gl_n(F))$ and $f\in C_c^\infty(\fh_n(F))$ match then so do $\cF f'$ and $\cF f$.
\end{theo}
$\blacksquare$

\section{End of the proof}

Set $\Phi=O(.,\bgln)-O(.,\bun)\in Orb(\gl_n)+Orb(\fh_n)$. Clearly, Theorem \ref{theo 1} is equivalent to the identical vanishing of $\Phi$. As $\Phi$ is a locally constant function on $\cA^{\rs}(F)$ it suffices to show that $\Phi(a)=0$ for the dense subset of $a\in \cA^{\rs}(F)$ such that $q(a)\neq 0$. By Lemma \ref{lem 1} we at least have
$$\displaystyle \Phi(a)=0 \mbox{ for every }a\in \cA^{\rs}(F) \mbox{ such that } \lvert q(a)\rvert\geqslant 1.$$
Hence, by \eqref{eqC} and since $\psi$ is unramified, we have
$$\displaystyle W\begin{pmatrix} 1 & t \\ 0 & 1\end{pmatrix}\Phi=\Phi$$
for every $t\in \mathfrak{p}_F^{-1}$. On the other hand, since $\psi$ is unramified we have $\cF \bgln=\bgln$ and $\cF\bun=\bun$ so that
$$\displaystyle W\begin{pmatrix} 0 & -1 \\ 1 & 0 \end{pmatrix}\Phi=O(.,\cF\bgln)-O(.,\cF\bun)=\Phi.$$
As $SL_2(F)$ is generated by $\begin{pmatrix} 0 & -1 \\ 1 & 0 \end{pmatrix}$ and $\begin{pmatrix} 1 & \mathfrak{p}_F^{-1}\\ 0 & 1 \end{pmatrix}$, it follows that $\Phi$ is fixed by the Weil representation. In particular
$$\displaystyle \psi(tq(a))\Phi(a)=W\begin{pmatrix} 1 & t \\ 0 & 1\end{pmatrix}\Phi(a)=\Phi(a)$$
for every $a\in \cA^{\rs}(F)$ and $t\in F$. When $q(a)\neq 0$, we can always find $t$ such that $\psi(tq(a))\neq 1$ hence $\Phi(a)=0$ and this ends the proof of Theorem \ref{theo 1}.

\section*{Acknowledgment}
I am grateful to Jean-Loup Waldspurger for a very careful reading and for correcting many inaccuracies in a first version of this paper. I also thank Hang Xue for useful comments and the referees for their thorough readings and many suggestions that have improved the level of exposition of this paper.


\begin{thebibliography}{9}
\bibitem[Bump]{Bump}
D. Bump: \textit{Automorphic forms and representations}, Cambridge Studies in Advanced Mathematics, 55. Cambridge University Press, Cambridge, 1997. xiv+574 pp. ISBN: 0-521-55098-X

\bibitem[Chau]{Chau}
P.-H. Chaudouard: \textit{On relative trace formulae: the case of Jacquet-Rallis}, prepublication 2017 available at \url{https://webusers.imj-prg.fr/~pierre-henri.chaudouard/Hanoi.pdf}

\bibitem[JL]{JL}
H. Jacquet, R.P. Langlands: \textit{Automorphic forms on $GL(2)$}, Lecture Notes in Math., Vol. 114, Springer-Verlag, Berlin and New York, 1976.

\bibitem[JR]{JR}
H. Jacquet, S. Rallis: \textit{On the Gross-Prasad conjecture for unitary groups}, in ``On certain $L$-functions", 205-264, Clay Math. Proc., 13, Amer. Math. Soc., Providence, RI, 2011.

\bibitem[KV]{KV}
D. Kazhdan, Y. Varshavsky: \textit{On endoscopic transfer of Deligne Lusztig functions}, Duke Math. J. 161 (2012), no. 4, 675-732.

\bibitem[Wa1]{Wa1}
J.-L. Waldspurger: \textit{Une formule des traces locale pour les alg\`ebres de Lie $p$-adiques}, J. Reine Angew. Math. 465 (1995) 41-99.

\bibitem[Wa2]{Wa2}
J.-L. Waldspurger: \textit{Le lemme fondamental implique le transfert}, Compositio Math. 105 (1997), 153-236.

\bibitem[Xi]{Xi}
J. Xiao: \textit{Endoscopic transfer for unitary Lie algebras}, prepublication 2018 arXiv:1802.07624

\bibitem[Yu]{Yu}
Z. Yun: \textit{The fundamental lemma of Jacquet and Rallis}, With an appendix by Julia Gordon. Duke Math. J. 156 (2011), no. 2, 167-227.

\bibitem[Zha1]{Zha1}
W. Zhang: \textit{Fourier transform and the global Gan-Gross-Prasad conjecture for unitary groups}, Ann. of Math. (2), 180(3):971-1049, 2014.
\end{thebibliography}
\end{document}